\documentclass[reqno]{amsart}
\begin{document}

\title
[Improvement Taykov's lower bound in an inequation between $C$ and $L$ norms]
{Improvement Taykov's lower bound\\ in an inequation between $C$ and $L$ norms\\ for trigonometric polynomials}
\author{\sc D.V. Gorbachev}
\date{10/17/2002}
\thanks{This research was supported by the RFBR under grant no.~00-01-00644 and 02-01-06563.
\endgraf
Matem. zametki (Russia), 2002 (in publishing).}
\email{dvg@uic.tula.ru}
\urladdr{http://home.uic.tula.ru/$\sim$gd030473}
\begin{abstract}
We give new lower asymptotical estimate of constant
\[
C_n=\sup\biggl\{\frac{\|t_n\|_{C(\mathbb T)}}{\|t_n\|_{L(\mathbb T)}}:t_n\text{ are real trigonometric
polynomials},\ \operatorname{deg}t_n<n\biggr\}
\]
as $n\to\infty$. This estimate improves known bound of L.V.Taykov (1965).
\end{abstract}
\maketitle

Suppose $e(t)=\exp(2\pi it)$, $c(t)=\cos(2\pi t)$, $s(t)=\sin(2\pi t)$. We consider
an extremum problem about behavior of the best constant as $n\to\infty$ in the
inequation between $C$ and $L$ norms for real trigonometric polynomials
$t_n(x)=\sum_{\nu\in\mathbb Z,\,|\nu|<n}\widehat t_n(\nu)e(\nu x)$:
\[C_n=\sup\limits_{t_n}\frac{\|t_n\|_C}{\|t_n\|_L},\] where
\[\|t_n\|_C=\max\limits_{|x|\le1/2}|t_n(x)|,\qquad
\|t_n\|_L=\int_{-1/2}^{1/2}|t_n(x)|\,dx.\]

Exact value of the constant $C_n$ is unknown for all $n\ge3$ [1] ($C_1=1$,
$C_2=2{,}12532\ldots$). S.B.~Stechkin (see [1,~2]) proved that $C_n=c\cdot n+o(n)$
as $n\to\infty$ ($c>0$). L.V.~Taykov [2] improved this result. He found following
bounds for $C_n$:
\begin{equation}\label{1}
(1{,}07995\ldots)\,n+O(1)=\frac2{\int_0^\pi\sin t\,dt/t}\,n+O(1)\le C_n\le
\frac{4G}\pi\,n+O(1)=(1{,}16624\ldots)\,n+O(1),
\end{equation}
where $G=\sum_{k=0}^\infty(-1)^k(2k+1)^{-2}=0{,}91596\ldots$ is Katalan's
constant. The lower bound of this result is obtained for polynomials of Rogozinsky
[3]
\begin{equation}\label{2}
R_n(x)=\sum_{|\nu|<n}c(\tfrac\nu{4n})\,e(\nu x)=
s(\tfrac1{4n})c(nx)(c(x)-c(\tfrac1{4n}))^{-1}.
\end{equation}

In this work we show new polynomials for improvement the lower bound in \eqref{1}.
By $q_k=k/2-1/4$ ($k\in\mathbb Z$)~  denote zeros of function $c(t)$;
$1/2<z_1<1<z_2<3/2$, $$
\omega(z)=\frac{(z_1^2-z^2)(z_2^2-z^2)}{(q_1^2-z^2)(q_2^2-z^2)(q_3^2-z^2)},\qquad
\omega_k(z)=(q_k^2-z^2)\omega(z)\quad(k=1,2,3). $$

Consider an even function \[\varphi(x)=\sum_{k=1}^3\omega_k(q_k)\sigma_k(x)\] such
that $\sigma_k(x)=(\pi/q_k)s(q_k(1-|x|))$ for $|x|<1$ and $\sigma_k(x)=0$ for
$|x|\ge1$.

Fourier transform $\widehat\sigma_k(z)=\int_{\mathbb R}\sigma_k(x)e(-zx)\,dx$
equals to $c(z)(q_k^2-z^2)^{-1}$. Therefore Fourier transform of the function
$\varphi$ equals to
\[\widehat\varphi(z)=\sum_{k=1}^3\omega_k(q_k)\,\frac{c(z)}{q_k^2-z^2}=\omega(z)c(z).\]

Suppose \[F_n(x)=\sum_{|\nu|<n}\varphi(\tfrac\nu n)e(\nu x).\] Using Poisson
summation formula $\sum_{\nu\in\mathbb Z}\widehat f(\nu-x)=\sum_{\nu\in\mathbb
Z}f(\nu)e(\nu x)$ for the function $f(x)=\varphi(x/n)$ (since
$\widehat\varphi(z)=O(z^{-2})$ ($|z|\to\infty$), $\varphi(x)=0$ ($|x|\ge1$)) and
identity \[\sum_{\nu\in\mathbb Z}\frac1{a^2-(z+\nu)^2}=\frac{(\pi/a)s(a)}{c(z)-c(a)},\]
we get a following expression for polynomial $F_n$
\begin{equation}\label{3}
F_n(x)=n\sum_{\nu\in\mathbb Z}\widehat\varphi(nx+n\nu)=
c(nx)\sum_{k=1}^3\frac{(\pi/q_k)\omega_k(q_k)s(\tfrac{q_k}n)}{c(x)-c(\tfrac{q_k}n)}.
\end{equation}

\smallskip\textsc{Theorem.} {\itshape
$\|F_n\|_C/\|F_n\|_L\ge\rho(\alpha,\beta)\,n+O(n^{-1}),$ where $$
\rho(z_1,z_2)=\widehat\varphi(0)
\biggl[4\biggl(\int_0^{q_1}+\int_{q_2}^{q_3}-\int_{z_1}^{z_2}\,\biggr)\widehat\varphi(z)\,dz\biggr]^{-1}.
$$}

\smallskip The function $\rho(z_1,z_2)$ can be expressed by an integral sine. Using
necessary conditions of extreme, we maximize this function on $1/2<z_1<1<z_2<3/2$
and obtain new asymptotic lower bound for the constant~$C_n$.

\smallskip\textsc{Corollary.} {\itshape $C_n\ge(1.08176\ldots)\,n+O(n^{-1}),$
$z_1^*=0{,}72096\ldots,$ $z_2^*=1{,}23305\ldots$}

\smallskip\textsc{Remark 1.} {\itshape For $z_1=q_2$, $z_2=q_3$ polynomials $F_n$
coincide with polynomials of Rogozinsky \eqref{2} \textup(accurate to a coefficient
$4\pi$\textup) $F_n(x)=4\pi R_n(x)$ and $\rho(q_2,q_3)=2(\int_0^\pi\sin t\,dt/t)^{-1};$
this fact give Taykov's lower bound.}

\smallskip\textsc{Proof.} Since $\widehat\varphi(z)=O(z^{-2})$ as $|z|\to\infty$, it
follows that $n\sum_{\nu\in\mathbb
Z\setminus\{0\}}\widehat\varphi(nx+n\nu)=O(n^{-1})$ uniformly for $|x|\le1/2$.
Therefore from \eqref{3} we have $$
\|F_n\|_C=\max\limits_{|x|\le1/2}\,\biggl|n\sum_{\nu\in\mathbb
Z}\widehat\varphi(nx+n\nu)\biggr|= \max\limits_{|z|\le
n/2}|n\widehat\varphi(z)|+O(n^{-1})\ge n\widehat\varphi(0)+O(n^{-1}). $$

Now we estimate $\|F_n\|_L$. We have
\begin{multline*}
\|F_n\|_L=\int_{-1/2}^{1/2}\,\biggl|n\sum_{\nu\in\mathbb Z}\widehat\varphi(nx+n\nu)\biggr|\,dx=
\int_{-1/2}^{1/2}|n\widehat\varphi(nx)|\,dx+O(n^{-1})
=\\=
\int_{-n/2}^{n/2}|\widehat\varphi(z)|\,dz+O(n^{-1})\le
\int_{\mathbb R}|\widehat\varphi(z)|\,dz+O(n^{-1}).
\end{multline*}

Therefore $\|F_n\|_C/\|F_n\|_L\ge\dfrac{\widehat\varphi(0)}{\int_{\mathbb
R}|\widehat\varphi(z)|\,dz}\,n+O(n^{-1})$.

We calculate $\int_{\mathbb R}|\widehat\varphi(z)|\,dz$. Suppose
\[E_+=[0,z_1]\cup[z_2,q_4]\cup(\cup_{k=2}^\infty[q_{2k+1},q_{2k+2}]),\qquad
E_-=(z_1,z_2)\cup(\cup_{k=2}^\infty(q_{2k},q_{2k+1})).\] Then $\mathbb R_+=E_+\cup
E_-$ and $\widehat\varphi(z)\ge0$ for $z\in E_+$; $\widehat\varphi(z)<0$ for $z\in
E_-$. $\widehat\varphi$ is even function. Hence we obtain
\begin{multline*}
\int_{\mathbb R}|\widehat\varphi(z)|\,dz=2\int_{\mathbb
R_+}|\widehat\varphi(z)|\,dz=
2\biggl(\int_0^{z_1}+\int_{z_2}^{q_4}+\sum_{k=2}^\infty\int_{q_{2k+1}}^{q_{2k+2}}-
\int_{z_1}^{z_2}-\sum_{k=2}^\infty\int_{q_{2k}}^{q_{2k+1}}
\,\biggr)\widehat\varphi(z)\,dz
=\\=
2\biggl(\int_{\mathbb R_+}-2\biggl(\int_{z_1}^{z_2}+
\sum_{k=2}^\infty\int_{k-1/4}^{k+1/4}\,\biggr)\biggr)\widehat\varphi(z)\,dz=
\varphi(0)-4\int_{z_1}^{z_2}\widehat\varphi(z)\,dz-
4\sum_{k=2}^\infty\int_{-1/4}^{1/4}\widehat\varphi(z+k)\,dz=\\=
\varphi(0)-4\int_{z_1}^{z_2}\widehat\varphi(z)\,dz-
4\int_0^{1/4}\biggl(\,\sum_{k=-\infty}^\infty\widehat\varphi(z+k)-\widehat\varphi(z)-
\widehat\varphi(z-1)-\widehat\varphi(z+1)\biggr)\,dz=\\=
\varphi(0)-4\int_{z_1}^{z_2}\widehat\varphi(z)\,dz-
4\int_0^{1/4}\varphi(0)\,dz+4\biggl(\int_0^{1/4}+\int_{-1}^{-3/4}+\int_{1}^{5/4}\,\biggr)
\widehat\varphi(z)\,dz=\\=
4\biggl(\int_0^{q_1}+\int_{q_2}^{q_3}-\int_{z_1}^{z_2}\,\biggr)\widehat\varphi(z)\,dz.
\end{multline*}

In these calculations we use parity
$\sum_{k=-\infty}^\infty\widehat\varphi(z+k)=\varphi(0)$; it following from Poisson
summation formula $\sum_{\nu\in\mathbb Z}\widehat\varphi(z+\nu)=\sum_{\nu\in\mathbb
Z}\varphi(\nu)e(-\nu z)$ and $\varphi(\nu)=0$ for $|\nu|\ge1$. Theorem is proved.

\smallskip\textsc{Remark 2.} {\itshape It is possible to improve new lower bound of
the constant $c,$ but it demands large calculations \textup(these calculations are not
included in this short article\textup).}

\smallskip
All computer calculations are obtained by mathematical package
\verb"Maple".

\end{document}